\documentclass{amsart}
\usepackage[dvips]{graphicx}
\usepackage{amscd}
\usepackage{amsmath}
\usepackage{amsxtra}
\usepackage{amsfonts}
\usepackage{amssymb}
\newtheorem{theorem}{Theorem}[section]
\newtheorem{corollary}[theorem]{Corollary}

\newtheorem{proposition}[theorem]{Proposition}
\theoremstyle{definition}
\newtheorem{definition}[theorem]{Definition}
\newtheorem{remark}[theorem]{Remark}

\newtheorem{example}[theorem]{Example}
\theoremstyle{remark}

\renewcommand{\theclaim}{\textup{\theclaim}}

\newtheorem*{acknowledgements}{Acknowledgements}

\numberwithin{equation}{section}

\newlength{\extramargin}
\def\openone

{\mathchoice

{\hbox{\upshape \small1\kern-3.3pt\normalsize1}}

{\hbox{\upshape \small1\kern-3.3pt\normalsize1}}

{\hbox{\upshape \tiny1\kern-2.3pt\SMALL1}}

{\hbox{\upshape \Tiny1\kern-2pt\tiny1}}}

\makeatletter

\newbox\ipbox

\newcommand{\ip}[2]{\left\langle #1\,|\,#2\right\rangle}

\newcommand{\diracb}[1]{\left\langle #1\mathrel{\mathchoice

{\setbox\ipbox=\hbox{$\displaystyle \left\langle\mathstrut
#1\right.$}

\vrule height\ht\ipbox width0.25pt depth\dp\ipbox}

{\setbox\ipbox=\hbox{$\textstyle \left\langle\mathstrut
#1\right.$}

\vrule height\ht\ipbox width0.25pt depth\dp\ipbox}

{\setbox\ipbox=\hbox{$\scriptstyle \left\langle\mathstrut
#1\right.$}

\vrule height\ht\ipbox width0.25pt depth\dp\ipbox}

{\setbox\ipbox=\hbox{$\scriptscriptstyle \left\langle\mathstrut
#1\right.$}

\vrule height\ht\ipbox width0.25pt depth\dp\ipbox}

}\right. }

\newcommand{\dirack}[1]{\left. \mathrel{\mathchoice

{\setbox\ipbox=\hbox{$\displaystyle \left.\mathstrut
#1\right\rangle$}

\vrule height\ht\ipbox width0.25pt depth\dp\ipbox}

{\setbox\ipbox=\hbox{$\textstyle \left.\mathstrut
#1\right\rangle$}

\vrule height\ht\ipbox width0.25pt depth\dp\ipbox}

{\setbox\ipbox=\hbox{$\scriptstyle \left.\mathstrut
#1\right\rangle$}

\vrule height\ht\ipbox width0.25pt depth\dp\ipbox}

{\setbox\ipbox=\hbox{$\scriptscriptstyle \left.\mathstrut
#1\right\rangle$}

\vrule height\ht\ipbox width0.25pt depth\dp\ipbox}

} #1\right\rangle}

\newcommand{\lonet}{L^{1}\left(  \mathbb{T}\right)}

\newcommand{\ltwot}{L^{2}\left(  \mathbb{T}\right)}

\newcommand{\linft}{L^{\infty}\left(  \mathbb{T}\right)}

\newcommand{\ltwor}{L^{2}\left(\mathbb{R}\right)}

\newcommand{\linfr}{L^{\infty}\left(\mathbb{R}\right)}

\input cyracc.def

\begin{document}
\title[Positive definite maps]{Positive definite maps, representations and frames}
\author{Dorin Ervin Dutkay}
\address{Department of Mathematics\\
The University of Iowa\\
14 MacLean Hall\\
Iowa City, IA 52242-1419\\
U.S.A.} \email{ddutkay@math.uiowa.edu}
\thanks{}
 \keywords{}

\begin{abstract}We present a unitary approach to the construction of representations and intertwining operators. We apply it to the $C^*$-algebras,
groups, Gabor type unitary systems and wavelets. We give an
application of our method to the theory of frames, and we prove a
general dilation theorem which is in turn applied to specific
cases, and we obtain in this way a dilation theorem for wavelets.

\end{abstract}\maketitle
\tableofcontents
\section{\label{Intro}Introduction}
 Engineering problems in time-frequency analysis of coherent vector expansions, Gabor bases, wavelets based on
scaling and integral translations, and multiresolution algorithms
in signal processing are generally not thought to be related to
operator algebras. In this paper, we show nonetheless that a
fundamental idea of Kolmogorov adds clarity to known constructions
in operator algebra theory, and moreover is the key to an
extension of recent results in the more applied areas that we
enumerated above. Of our original results (see section 5 below) we
highlight a new algorithm for the construction of certain
orthonormal frames of wavelet type. Our paper proposes a general
method of construction of representations of various algebraic
structures as operators on Hilbert spaces. Our goal is to show how
some well known constructions of representations fit into the same
framework and are consequences of a general result. Among the
structures considered, we mention $C^*$-algebras, groups, Gabor
type unitary systems and wavelet representations.
\par
In operator theory, the GNS construction producing representations
of $C^*$-algebras is a fundamental tool (see \cite{BraRo}). In
harmonic analysis unitary representations of groups can be
constructed when a function of positive type is present (see
\cite{Fol}). Representations are ubiquitous also in the theory of
wavelets and frames (see \cite{HL}, \cite{Jor}). We will see how
these various results have in fact a common ground - a classical
theorem of Kolmogorov (theorem \ref{th1_2}), also known in the
literature as the Kolmogorov decomposition of positive definite
kernels. We follow here the ideas introduced in \cite{EvLe}. It is
shown there that the Kolmogorov theorem gives a unified treatment
of several important dilation theorems such as the GNS-Stinespring
construction for $C^*$-algebras, the Naimark-Sz. Nagy unitary
dilation of positive definite functions on groups, the
construction of Fock spaces and the algebras of canonical
commutation and anticommutation relations. Kolmogorov's result was
used also by Sz. Nagy and C. Foias in dilation theory, for the
commutant lifting theorem (\cite{SzF68}, \cite{SzF70}) which in
turn was a key idea used by D. Sarason to obtain a solution to the
Nevanlinna-Pick interpolation problem (\cite{Sar67}). For a more
complete account of the history and applications of Komogorov's
result, we refer to \cite{C96}.
\par
We will indicate how this technique can be used also for
construction of wavelet representations and Gabor type unitary
systems.
\par
More general constructions for Hermitian kernels are also possible
and they are based on Krein spaces (see \cite{C97}).
\par
In section \ref{Positive} we review the general result of
Kolmogorov and we show how it can be used for the GNS construction
and for positive definite maps on groups.
 Then we apply it
to Gabor type unitary systems and we obtain unitary
representations and for wavelets we get the cyclic representations
introduced in \cite{Jor}.
\par
Section \ref{Intertwining} concerns operators compatible with the
representations defined in section \ref{Positive}, called
intertwining operators. Again, the starting point is a general
theorem (theorem \ref{th3_2}). We consider some particular cases
and study how the intertwining operators will be compatible with
the additional structure that appears.
\par
In section \ref{Frames} we analyze some connections between
representations and frames. We recall that a set $\{x_n\,|\,
n\in\mathbb{N}\}$ of vectors in a Hilbert space $H$ is called a
frame for the Hilbert space $H$ if there are some positive
constants $A$ and $B$ such that
$$A\|f\|^2\leq\sum_{n\in\mathbb{N}}|\ip{f}{x_n}|^2\leq B\|f\|^2,\quad(f\in H).$$
When $A=B=1$ the we call it normalized tight frame.
\par
It is known that any normalized tight frame is the projection of
an orthonormal basis of a bigger Hilbert space (see \cite{HL}). We
will prove that the normalized tight frames can be dilated to
orthonormal bases in a way that is compatible with the
representations defined in section \ref{Positive}. We will get as
immediate consequences the dilation theorems for groups and Gabor
type unitary systems introduced in \cite{HL}.
\par
In the last section we consider the case of wavelets obtained from
a multiresolution analysis. It is known (see \cite{Dau92}) that,
unless some restrictions are imposed on the low-pass filter that
starts the MRA construction, the wavelets obtained do not form an
orthonormal basis but a normalized tight frame. Since such frames
can be dilated to orthonormal bases, a natural question would be
if the dilation preserves the multiresolution structure. The
answer is affirmative and it is given in theorem \ref{th5_2} and
more concretely in theorem \ref{th5_3}. In this way we obtain
"wavelets" in a Hilbert space bigger then  $\ltwor$.

\section{\label{Positive}Positive definite maps and representations}
We begin this section with a general result of Kolmogorov
(\cite{EvLe}, \cite{EvKa}). Then we consider several structures
and show how to obtain representations from this general theorem.
\begin{definition}\label{def1_1}
Let $X$ be a nonempty set. We say that a map $K:X\times
X\rightarrow\mathbb{C}$ is positive definite, and denote this by
$0\leq K$, if
$$\sum_{i,j=1}^nK(x_i,x_j)\xi_i\overline{\xi}_j\geq0,\quad(n\in\mathbb{N}, x_i\in X, \xi_i\in\mathbb{C}\mbox{ for all }i\in\{1,...,n\}).$$
\end{definition}

\begin{theorem}\label{th1_2}{\bf [Kolmogorov's theorem]} If $K:X\times X\rightarrow\mathbb{C}$ is positive definite then there
exists a Hilbert space $H_K$ and a map $v_K:X\rightarrow H_K$ such
that the linear span of $\{v_K(x)\,|\, x\in X\}$ is dense in $H_K$
and
$$\ip{v_K(x)}{v_K(y)}=K(x,y),\quad(x,y\in X).$$
Moreover, $H_K$ and $v_K$ are unique up to unitary isomorphisms.
\end{theorem}
\begin{remark} Kolmogorov's theorem is valid also for
operator-valued positive definite maps and in this form it can be
applied for the Stinespring construction and the Naimark-Sz.Nagy
dilation. For details consult \cite{EvLe} and \cite{EvKa}.
\par
In this paper, for the application to wavelets and Gabor frames,
we will need only the more particular version of Kolmogorov's
theorem that we mentioned before.
\end{remark}
\begin{definition}\label{def1_3}
If $K:X\times X\rightarrow\mathbb{C}$ is positive definite then we
call $[H_K,v_K]$ the representation associated to $K$.
\end{definition}
\par
We note that Kolmogorov's theorem is purely set theoretic; there
is no structure on $X$. We expect that, if $X$ has some additional
structure on it and if we assume some compatibility between the
positive definite map $K$ and this structure, then the
representation associated to $K$ will also be in agreement with
the structure of $X$. In the next examples we will see that this
is indeed the case and we review the technique in the case of
$C^*$-algebras and groups.

\begin{example}\label{ex1_4}{\bf [$C^*$-algebras and the GNS construction]}
We consider now the case when $X=\mathcal{A}$ is a $C^*$-algebra
and prove that we can obtain the well known GNS construction from
Kolmogorov's theorem.
\end{example}

\begin{theorem}\label{th1_5}{\bf [The GNS construction]}
If $\mathcal{A}$ is a $C^*$-algebra and $\varphi$ is a positive
linear functional on $\mathcal{A}$, then there exists a
representation $\pi$ of $\mathcal{A}$ on a Hilbert space $H$, that
has a cyclic vector $\xi_0\in H$ such that
$$\ip{\pi(x)\xi_0}{\xi_0}=\varphi(x),\quad(x\in\mathcal{A}).$$
\end{theorem}
\begin{proof}
The idea is to define
$K:\mathcal{A}\times\mathcal{A}\rightarrow\mathbb{C}$ by
$$K(x,y)=\varphi(y^*x),\quad(x,y\in\mathcal{A}).$$
 We can use Kolmogorov's theorem to obtain the Hilbert
space $H_K$ and the map $v_K:\mathcal{A}\rightarrow H_K$.
\par
For a fixed $x\in\mathcal{A}$, define the operator $\pi(x)$ as
follows:
$$\pi(x)(v_K(y))=v_K(xy),\quad(y\in\mathcal{A}),$$
and extend by linearity. Then everything checks out.
\end{proof}

\begin{example}\label{ex1_6}{\bf [Groups and unitary representations]}
Take $X=G$ a group. We call $K:G\times G\rightarrow\mathbb{C}$ a
group positive definite map if $0\leq K$ and
$$K(x,y)=K(zx,zy),\quad(x,y,z\in G).$$
We note that such a positive definite map $K$ is uniquely
determined by its restriction $\phi(x)=K(x,1)$ and $\phi$ is a
function of positive type (see \cite{Fol}). The proof of theorem
\ref{th1_7} will show how the well known correspondence between
functions of positive type and unitary representations of groups
can be regarded as a consequence of Kolmogorov's theorem.
\end{example}

\begin{theorem}\label{th1_7}
Let $G$ be a group and $K$ a group positive definite map on $G$.
Then there exists a unitary representation $\pi_K$ of $G$ on a
Hilbert space $H_K$ with a cyclic vector $\xi_0\in H_K$ such that
$$\ip{\pi_K(x)\xi_0}{\pi_K(y)\xi_0}=K(x,y),\quad(x,y\in G).$$
\end{theorem}

\begin{proof}
The proof works exactly as in the case of $C^*$-algebras: consider
$[H_K,v_K]$ the representation associated to $K$ by Kolmogorov's
theorem. Define the operators $\pi_K(x)$ for $x\in G$ as follows:
$$\pi_K(x)(v_K(y))=v_K(xy),\quad(x,y\in G)$$
and extend by linearity.
\end{proof}

\begin{remark}\label{rem1_8}
Note that in the proof of theorem \ref{th1_7} we used the
representation associated to $K$ and we see that, when $K$ is a
group positive definite map, this representation has the unitary
representation $\pi_K$ attached to it. The same observation can be
done for the GNS construction: the representation of the
$C^*$-algebra is attached to the representation $v_K$. This
confirms our expectation: when the positive definite map has some
compatibility with the existent structure on $X$, this
compatibility projects a nice structure on the associated
representation $[H_K,v_K]$. This is the idea that we use
throughout this section.
\end{remark}

\begin{example}\label{ex1_9}{\bf [Gabor type unitary systems]}
We recall that a Gabor system is associated to two positive
constants $a,b>0$ and a function $g\in\ltwor$ and is defined by
$$g_{m,n}(\xi)=e^{2\pi imb\xi}g(\xi-na),\quad(\xi\in\mathbb{R}).$$
The Gabor systems are one of the major subjects in the study of
frames and wavelet theory. If we define the unitary operators
$U,V$ on $\ltwor$,
$$(Uf)(\xi)=e^{2\pi ib}f(\xi),\quad(f\in\ltwor),$$
$$(Vf)(\xi)=f(\xi-a),\quad(f\in\ltwor),$$
then $g_{m,n}=U^mV^ng$, $(m,n\in\mathbb{Z})$, and $U$ and $V$
satisfy the relation
$$UV=e^{2\pi iab}VU.$$
Following \cite{HL}, if $U$ and $V$ are unitary operators on a
Hilbert space $H$ that verify the relation
$$UV=\lambda VU$$
for some unimodular scalar $\lambda$, we then call $\{U^mV^n\,|\,
m,n\in\mathbb{Z}\}$ a Gabor type unitary system. We will prove
that these systems fit into our general framework and we construct
representations for them.
\end{example}

\begin{theorem}\label{th1_10}
Suppose $\lambda$ is a unimodular scalar and
$K:\mathbb{Z}^2\times\mathbb{Z}^2\rightarrow\mathbb{C}$ is a
positive definite map satisfying
\begin{equation}\label{eq1_10_1}
K((m+1,n),(m'+1,n'))=K((m,n),(m',n')),\quad(m,n,m',n'\in\mathbb{Z}),
\end{equation}
\begin{equation}\label{eq1_10_2}
K((m,n+1),(m',n'+1))=\lambda^{m-m'}K((m,n),(m',n')),\quad(m,n,m',n'\in\mathbb{Z}).
\end{equation}
Then, on the Hilbert space $H_K$, there are unitaries $U,V$ and a
vector $\xi_0\in H_K$ such that
\begin{eqnarray}
 &UV=\lambda VU& \label{eq1_10_3}\\
 &\{U^mV^n\xi_0\,|\,m,n\in\mathbb{Z}\}\mbox{ is dense in }H_K& \label{eq1_10_4}\\
 &\ip{U^mV^n\xi_0}{U^{m'}V^{n'}\xi_0}=K((m,n),(m',n')),\quad(m,n,m',n'\in\mathbb{Z}).& \label{eq1_10_5}
\end{eqnarray}
Moreover, this representation is unique up to unitary isomorphism.
\end{theorem}

\begin{proof}
Let $v_K:\mathbb{Z}^2\rightarrow H_K$ be the representation
associated to $K$. Define the operators $U$ and $V$ as follows:
$$U(v_K(m,n))=v_K(m+1,n),\quad(m,n\in\mathbb{Z}),$$
$$V(v_K(m,n))=\lambda^{-m}v_K(m,n+1),\quad(m,n\in\mathbb{Z}),$$
and then extend by linearity. We check that $U,V$ are well defined
and isometric. Take $a_i\in\mathbb{C}$,
$(m_i,n_i)\in\mathbb{Z}^2$, $(i\in\{1,...,p\})$.
$$\ip{V(\sum_{i=1}^pa_iv_K(m_i,n_i))}{V(\sum_{i=1}^pa_iv_K(m_i,n_i))}=$$
$$=\ip{\sum_{i=1}^pa_i\lambda^{-m_i}v_K(m_i,n_i+1)}{\sum_{i=1}^pa_i\lambda^{-m_i}v_K(m_i,n_i+1)}$$
$$=\sum_{i,j=1}^pa_i\overline{a}_j\lambda^{-m_i}\overline{\lambda^{-m_j}}K((m_i,n_i+1),(m_j,n_j+1))$$
$$=\sum_{i,j=1}^pa_i\overline{a}_jK((m_i,n_i),(m_j,n_j))
=\ip{\sum_{i=1}^pa_iv_K(m_i,n_i)}{\sum_{i=1}^pa_iv_K(m_i,n_i)}.$$
A similar calculation shows that $U$ is well defined and
isometric. Since the linear span of the vectors $v_K(m,n)$ is
dense in $H_K$, we can extend $U$ and $V$ to unitaries on $H_K$.
\par
Next, we check (\ref{eq1_10_3}).Take $(m,n)\in\mathbb{Z}^2$.
\begin{align*}
UVv_K(m,n)&=U(\lambda^{-m}v_K(m,n+1))=\lambda^{-m}v_K(m+1,n+1)\\
&=\lambda V(v_K(m+1,n))=\lambda VU(v_K(m,n)),
\end{align*}
and (\ref{eq1_10_3}) follows by density. Also, note that, if
$\xi_0=v_K(0,0)$, then
$$U^mV^n\xi_0=U^mV^nv_k(0,0)=U^mv_K(0,n)=v_K(m,n),\quad(m,n\in\mathbb{Z}).$$
This will imply (\ref{eq1_10_4}) and (\ref{eq1_10_5}).The
uniqueness is a consequence of the uniqueness part of Kolmogorov's
theorem.
\end{proof}

\begin{remark}\label{rem1_11}
Any Gabor type unitary system $U$, $V$ on a Hilbert space $H$,
that has a vector $\xi_0\in H$ with the property that the linear
span of
$$\{U^mV^n\xi_0\,|\,m,n\in\mathbb{Z}\}$$
is dense in $H$, gives rise to a positive definite map $K$ on
$\mathbb{Z}^2$ that satisfies (\ref{eq1_10_1}),(\ref{eq1_10_2}) as
follows:
$$K((m,n),(m',n'))=\ip{U^mV^n\xi_0}{U^{m'}V^{n'}\xi_0},\quad(m,n,m',n'\in\mathbb{Z}).$$
(\ref{eq1_10_1}),(\ref{eq1_10_2}) are just immediate consequences
of the fact that $U$ and $V$ are unitary and $UV=\lambda VU$.
\end{remark}

\begin{example}\label{ex1_12}{\bf [Wavelet representations]}
We recall briefly some facts about wavelet representations.
Wavelet theory deals with two unitary operators $U$ and $T$ on
$\ltwor$, corresponding to the integer $N\geq 2$ called the scale:
$$Uf(x)=\frac{1}{\sqrt{N}}f\left(\frac{x}{N}\right),\, Tf(x)=f(x-1),\quad(x\in\mathbb{R},f\in\ltwor).$$
 A wavelet is a function $\psi\in\ltwor$ such that
$$\{U^mT^n\psi\,|\,m,n\in\mathbb{Z}\}$$
is an orthonormal basis for $\ltwor$. One way to construct
wavelets is by multiresolutions and scaling functions (see
\cite{Dau92}). Scaling functions satisfy equations of the form
\begin{equation}\label{eq1_12_1}
U\varphi=\sum_{k\in\mathbb{Z}}a_KT^k\varphi,
\end{equation}
where $a_k$ are complex coefficients.
\par
The scaling equation can be reformulated using representations.
There is a representation of $\linft$ ($\mathbb{T}$ is the unit
circle) on $\ltwor$ given by
$$\widehat{\left(\pi(f)\xi\right)}=f\widehat{\xi},\quad(f\in\linfr,\xi\in\ltwor)$$
($\widehat{\xi}$ denotes the Fourier transform of $\xi$ and
functions on $\mathbb{T}$ are identified with $2\pi$-periodic
functions on $\mathbb{R}$). Using this representation,
(\ref{eq1_12_1}) can be rewritten as
$$U\varphi=\pi(m_0)\varphi,$$
$m_0(e^{-i\theta})=\sum_{k\in\mathbb{Z}}a_ke^{-ik\theta}$ is
called a low-pass filter. Also, the representation satisfies
$$U\pi(f)U^{-1}=\pi\left(f(z^N)\right),\quad(f\in\linft).$$
$(U,\pi,\ltwor,\varphi)$ is called the wavelet representation with
scaling function $\varphi$.
\par
The wavelet theory has shown a strong interconnection between
properties of the scaling function $\varphi$ and spectral
properties of the transfer operator associated to the low-pass
filter $m_0$:
$$R_{m_0,m_0}f(z)=\frac{1}{N}\sum_{w^N=z}|m_0|^2(w)f(w),\quad(z\in\mathbb{T},f\in\lonet).$$
where $\mathbb{T}$ is endowed with the normalized Haar measure.
For more information on this we refer the reader to \cite{BraJo}.
In particular, functions that are harmonic with respect to
$R_{m_0,m_0}$, i.e. $R_{m_0,m_0}h=h$, play an important role in
the theory.
\par
We recall here a theorem from \cite{Jor} which establishes the
link between functions which are harmonic with respect to
$R_{m_0,m_0}$ and wavelet representations, because it is another
particularized instance of Kolmogorov's theorem.
\end{example}

\begin{theorem}\label{th1_13}
If $m_0\in\linft$ is non-singular (i.e. it doesn't vanish on a set
of positive measure) and $h\in\lonet$, satisfies
$$R_{m_0,m_0}h=h,\, h\geq0,$$
then there exists a Hilbert space $H_h$, a representation $\pi_h$
of $\linft$ on $H_h$, a unitary $U_h$ on $H_h$ and a vector
$\varphi_h\in H_h$ such that
$$\overline{\mbox{span}}\{U_h^{-n}\pi_h(f)\varphi_{h}\,|\,n\in\mathbb{N},f\in\linft\}=H_h;$$
$$U_h\pi_h(f)U_h^{-1}=\pi_h(f(z^N)),\quad(f\in\linft);$$
$$U_h\varphi_h=\pi_h(m_0)\varphi_h;$$
$$\ip{\pi_h(f)\varphi_h}{\varphi_h}=\int_{\mathbb{T}}fh\,d\mu.$$
Moreover, this is unique up to unitary equivalence.
\end{theorem}

\begin{proof}
We give here only a sketch of the proof that uses Kolmogorov's
theorem, the rest are calculations wich can be found in
\cite{Jor}.
\par
Let
$$X=\{(f,n)\,|\, f\in\linft,n\in\mathbb{N}\}.$$
We want to define a positive definite map $K$ on $X$ such that in
the end $v_K(f,n)=U_h^{-n}\pi_h(f)\varphi_h$. Then, we must have
\begin{align*}
K((f,n),(g,m))&=\ip{U_h^{-n}\pi_h(f)\varphi_h}{U_h^{-m}\pi_h(g)\varphi_h}\\
&=\ip{U_h^{m}\pi_h(f)\varphi_h}{U_h^{n}\pi_h(g)\varphi_h}\\
&=\ip{\pi_h(f(z^{N^m})m_0^{(m)}(z))\varphi_h}{\pi_h(g(z^{N^n})m_0^{(n)}(z))\varphi_h}\\
&=\int_{\mathbb{T}}f(z^{N^m})m_0^{(m)}(z)\overline{g}(z^{N^n})\overline{m_0^{(n)}}(z)h\,d\mu,
\end{align*}
where $m_0^{(m)}(z)=m_0(z)m_0(z^N)...m_0(z^{N^{m-1}})$.
\par
So we have to define, for $(f,n),(g,m)\in X$,
$$K((f,n),(g,m))=\int_{\mathbb{T}}f(z^{N^m})m_0^{(m)}(z)\overline{g}(z^{N^n})\overline{m_0^{(n)}}(z)h\,d\mu.$$
$K$ can be checked to be positive definite so it induces a
representation $(H_h,v_h)$, according to Kolmogorov's thoerem.
Then, define $\varphi_h=v_h(1,0)$,
$$U_hv_h(f,0)=v_h(f(z^N)m_0,0),\quad(f\in\linft),$$
$$U_hv_h(f,n)=(f,n-1),\quad(n\geq1,f\in\linft),$$
and extend by linearity and density.
$$\pi_h(f)v_h(g,n)=(f(z^{N^n})g(z),n),\quad(f,g\in\linft,n\in\mathbb{N}),$$
and extend by linearity and density.
\par
Everything can be checked out as the reader may see in \cite{Jor}.
\end{proof}

\section{\label{Intertwining}Intertwining operators}
In the previous section we saw how positive definite maps induce
representations on Hilbert spaces. Now we will show that
intertwining operators can be constructed in a similar way from
maps $L:X\times X\rightarrow\mathbb{C}$ which satisfy some
boundedness condition. We will also see that, when $X$ has some
structure on it and $L$ is compatible with this structure, then
the intertwining operator induced by $L$ will be compatible with
the extra structure existent on the induced representations, i.e.
the operator is indeed intertwining.
\par
The format of this section is similar to the format of the
previous one. We begin with a general, set theoretic result and
then particularize it to various structures to obtain more
information.

\begin{definition}\label{def3_1}
Consider two positive definite maps $K,K':X\times
X\rightarrow\mathbb{C}$ and $L:X\times X\rightarrow\mathbb{C}$
(not necessarily positive definite). We say that $L$ is bounded
with respect to $K$ and $K'$ if there is a constant $c>0$ such
that
\begin{equation}\label{eq3_1}
\left|\sum_{i=1}^m\sum_{j=1}^nL(x_i,y_j)\xi_i\overline{\eta}_j\right|^2\leq
c\left(\sum_{i,i'=1}^mK(x_i,x_{i'})\xi_i\overline{\xi}_{i'}\right)
\left(\sum_{j,j'=1}^mK'(y_j,y_{j'})\eta_j\overline{\eta}_{j'}\right)
\end{equation}
for all $x_i,y_j\in X$, $\xi_i,\eta_j\in\mathbb{C}$,
$i\in\{1,...,m\},j\in\{1,...,n\}$. We denote this by
$$L^2\leq cKK'.$$
\end{definition}

\begin{theorem}\label{th3_2}
Suppose $X$ is a nonempty set and $K,K'$ are positive definite
maps on $X$. If $L:X\times X\rightarrow\mathbb{C}$ and $L^2\leq
cKK'$ for some $c>0$, then there exists a unique bounded linear
operator $S:H_K\rightarrow H_{K'}$ such that
\begin{equation}\label{eq3_2_1}
\ip{Sv_K(x)}{v_{K'}(y)}=L(x,y),\quad(x,y\in X).
\end{equation}
($(H_K,v_K)$, $(H_{K'},v_{K'})$ are the representation induced by
$K$ and $K'$ respectively, according to Kolmogorov's theorem).
Moreover, $\|S\|\leq\sqrt{c}$. Conversely, if $S:H_K\rightarrow
H_{K'}$ is a bounded linear operator, then there is a unique map
$L:X\times X\rightarrow\mathbb{C}$ with $L^2\leq\|S\|^2KK'$ that
satisfies (\ref{eq3_2_1}).
\end{theorem}

\begin{proof}
Define $B:H_K\times H_{K'}\rightarrow\mathbb{C}$ as follows: for
$x_i,y_j\in X,\xi_i,\eta_j\in\mathbb{C})$,
$$B(\sum_{i=1}^n\xi_iv_K(x_i),\sum_{j=1}^n\eta_jv_{K'}(y_j))=\sum_{i=1}^m\sum_{j=1}^n\xi_i\overline{\eta}_jL(x_i,y_j).$$
Because $L^2\leq cKK'$, we have
$$\left|B(\sum_{i=1}^n\xi_iv_K(x_i),\sum_{j=1}^n\eta_jv_{K'}(y_j))\right|^2\leq c\left(\sum_{i,i'=1}^mK(x_i,x_{i'})\xi_i\overline{\xi}_{i'}\right)\cdot$$
$$\cdot\left(\sum_{j,j'=1}^mK'(y_j,y_{j'})\eta_j\overline{\eta}_{j'}\right)=c\left\|\sum_{i=1}^n\xi_iv_K(x_i)\right\|_{H_K}^2\left\|\sum_{j=1}^n\eta_jv_{K'}(y_j)\right\|_{H_{K'}}^2$$
This shows that $B$ is a well defined bounded sesquilinear map
which can be extended (by the density properties of $v_K$ and
$V_{K'}$) to a bounded sesquilinear map $B:H_K\times
H_{K'}\rightarrow\mathbb{C}$. Then there exists a bounded linear
operator $S:H_K\rightarrow H_{K'}$ such that $\|S\|\leq\sqrt{c}$
and
$$B(v_1,v_2)=\ip{Sv_1}{v_2},\quad(v_1\in H_K,v_2\in H_{K'}).$$
In particular, one obtains (\ref{eq3_2_1}).
\par
The uniqueness is clear because the spans of $\{v_K(x)\,|\, x\in
X\}$ and $\{v_{K'}(y)\,|\, y\in X\}$ are dense. The converse is
also easy, one needs to check that the map $L$ defined by
(\ref{eq3_2_1}) satisfies $L^2\leq \|S\|^2KK'$, but this is a
consequence of Schwarz's inequality.
\end{proof}

\begin{definition}\label{def3_2_1}
We call the operator $S$ associated to $L$ in theorem \ref{th3_2},
the intertwining operator associated to $L$.
\end{definition}

\par
We will also be interested in subrepresentations and in the
commutant of a representation. In these instances we will work
with only one positive definite map $K$. We give here a definition
which will be appropriate for these situations.

\begin{definition}\label{def3_3}
Consider $K,K'$, two positive definite maps on a nonempty set $X$
and a constant $c>0$. We denote
$$K'\leq cK$$
if, for all $x_i\in X$ and $\xi_i\in\mathbb{C}$,
$(i\in\{1,...,n\}$),
$$\sum_{i,j=1}^nK'(x_i,x_j)\xi_i\overline{\xi}_j\leq c\sum_{i,j=1}^nK(x_i,x_j)\xi_i\overline{\xi}_j.$$
\end{definition}

\begin{proposition}\label{prop3_4}
If $K$ and $K'$ are positive definite maps and $c>0$ then $K'\leq
cK$ if and only if $K'^2\leq c^2KK$.
\end{proposition}
\begin{proof}
Suppose $K'\leq cK$. Take $x_i,y_j\in X$,
$\xi_i,\eta_j\in\mathbb{C}$.
$$\left|\sum_{i=1}^m\sum_{j=1}^nK'(x_i,y_j)\xi_i\overline{\eta}_j\right|^2=\left|\ip{\sum_{i=1}^m\xi_iv_{K'}(x_i)}{\sum_{j=1}^n\eta_jv_{K'}(y_j)}\right|^2$$
$$\leq\left\|\sum_{i=1}^m\xi_iv_{K'}(x_i)\right\|_{H_{K'}}^2\left\|\sum_{j=1}^n\eta_jv_{K'}(y_j)\right\|_{H_{K'}}^2$$
$$=\left(\sum_{i,i'=1}^mK'(x_i,x_{i'})\xi_i\overline{\xi}_{i'}\right)
\left(\sum_{j,j'=1}^mK'(y_j,y_{j'})\eta_j\overline{\eta}_{j'}\right)$$
$$\leq c^2\left(\sum_{i,i'=1}^mK(x_i,x_{i'})\xi_i\overline{\xi}_{i'}\right)
\left(\sum_{j,j'=1}^mK(y_j,y_{j'})\eta_j\overline{\eta}_{j'}\right).$$
Hence $K'^2\leq c^2 KK$. Conversely, if $K'^2\leq c^2 KK$ then
just take $m=n$, $x_i=y_i$, $\xi_i=\eta_i$ in (\ref{eq3_1}) to
obtain exactly $K'\leq cK$.
\end{proof}

\begin{corollary}\label{cor3_5}
Suppose $K$ is positive definite on $X$. Then, for every positive
definite map $K'$ with $K'\leq cK$ for some $c>0$, there exists a
unique positive operator $S:H_K\rightarrow H_K$ with
\begin{equation}\label{eq3_2}
\ip{Sv_K(x)}{v_K(y)}=K'(x,y),\quad(x,y\in X).
\end{equation}
Moreover, $\|S\|\leq c$. Conversely, for every positive operator
$S$ on $H_K$ there is a unique positive definite map on $X$ that
satisfies (\ref{eq3_2}). In addition, $K'\leq \|S\|K$.
\end{corollary}
\begin{proof}
Using proposition \ref{prop3_4} and theorem \ref{th3_2}, we find
an operator $S$ on $H_K$ that satisfies (\ref{eq3_2}) and
$\|S\|\leq c$. $S$ is positive because
$$\ip{S(\sum_{i=1}^n\xi_iv_K(x_i))}{\sum_{i=1}^n\xi_iv_K(x_i)}=\sum_{i,j=1}^n\xi_i\overline{\xi}_jK'(x_i,x_j)\geq 0$$
and $\{v_K(x)\,|\,x\in X\}$ span a dense subspace of $H_K$.
\par
For the converse, when $S$ is given, theorem \ref{th3_2} shows
that there is a $K'$ satisfying (\ref{eq3_2}) and
$K'^2\leq\|S\|^2KK$. $K'$ is positive because $S$ is, and
proposition \ref{prop3_4} implies $K'\leq\|S\|K$.
\end{proof}
In the remainder of this section we apply theorem \ref{th3_2} to
the situations when $X$ has some additional structure on it and
see how the intertwining operators are in compliance with the
extra structure of the representations.

\begin{example}\label{ex3_6}{\bf [$C^*$-algebras]}
Consider now $X=\mathcal{A}$, a $C^*$-algebra. We saw in example
\ref{ex1_4} that, when the positive definite map
$K:\mathcal{A}\times\mathcal{A}\rightarrow\mathbb{C}$ is given by
a positive functional $\varphi:\mathcal{A}\rightarrow\mathbb{C}$,
$$K(x,y)=\varphi(y^*x),\quad(x,y\in\mathcal{A},$$
then the representation induced by $K$ has the GNS construction
attached to it. We want to see for what functions
$L:\mathcal{A}\times\mathcal{A}\rightarrow\mathbb{C}$ the
associated intertwining operator will intertwine the GNS
representations.
\end{example}

\begin{theorem}\label{th3_7}
Let $\mathcal{A}$ be a $C^*$-algebra and $\varphi,\varphi'$ two
positive functionals on $\mathcal{A}$. Suppose that
$\varphi_0:\mathcal{A}\rightarrow\mathbb{C}$ is linear and
$\varphi_0^2\leq c\varphi\varphi'$ for some $c>0$ in the sense
that
\begin{equation}\label{eq3_7_1}
|\varphi_0(y^*x)|^2\leq
c\varphi(x^*x)\varphi'(y^*y),\quad(x,y\in\mathcal{A}).
\end{equation}
Then there exists a unique bounded operator
$S:H_{\varphi}\rightarrow H_{\varphi'}$ such that
\begin{equation}\label{eq3_7_2}
S\pi_{\varphi}(x)=\pi_{\varphi'}(x)S,\quad(x\in\mathcal{A}),
\end{equation}
\begin{equation}\label{eq3_7_3}
\ip{S\pi_{\varphi}(x)\xi_0}{\pi_{\varphi'}(y)\xi_0'}=\varphi_0(y^*x),\quad(x,y\in\mathcal{A}).
\end{equation}
(Here $(H_{\varphi},\pi_{\varphi},\xi_0)$ and
$(H_{\varphi'},\pi_{\varphi'},\xi_0')$ are the GNS representations
associated to $\varphi$ and $\varphi'$ respectively (see theorem
\ref{th1_5}).) Moreover $\|S\|\leq\sqrt{c}$. Conversely, if
$S:H_{\varphi}\rightarrow H_{\varphi'}$ is a bounded operator that
satisfies (\ref{eq3_7_2}) then there is a unique linear map
$\varphi_0:\mathcal{A}\rightarrow\mathbb{C}$ that satisfies
(\ref{eq3_7_3}). In addition (\ref{eq3_7_1}) holds with
$c=\|S\|^2$.
\end{theorem}
\begin{proof}
Let
$K_{\varphi},K_{\varphi'}:\mathcal{A}\times\mathcal{A}\rightarrow\mathbb{C}$,
$$K_{\varphi}(x,y)=\varphi(y^*x),K_{\varphi'}(x,y)=\varphi'(y^*x),\quad(x,y\in\mathcal{A}).$$
Recall that
$H_{\varphi}=H_{K_{\varphi}},H_{\varphi'}=H_{K_{\varphi'}},\pi_{\varphi}(x)\xi_0=v_{K_{\varphi}}(x),
\pi_{\varphi'}(x)\xi_0=v_{K_{\varphi'}}(x)$ (see the proof of
theorem \ref{th1_5}).
\par
Define $L(x,y)=\varphi_0(y^*x)$ for $x,y\in\mathcal{A}$. Then
(\ref{eq3_7_1}) implies $L^2\leq cK_{\varphi}K_{\varphi'}$.
Theorem \ref{th3_2} gives an operator $S$ with
$$\ip{Sv_{K_{\varphi}}(x)}{v_{K_{\varphi'}}(y)}=L(x,y),\quad(x,y\in\mathcal{A});$$
then one checks that $S$ satisfies all the requirements.
\end{proof}
As a corrolary we deduce a basic fact about positive operators in
the commutant of the GNS representation (see\cite{BraRo}).
\begin{corollary}\label{cor3_8}
Let $\varphi,\varphi'$ be two positive functionals on a
$C^*$-algebra $\mathcal{A}$, $\varphi'\leq c\varphi$ for some
$c>0$ (i.e. $\varphi'(x)\leq c\varphi(x)$ for all positive
$x\in\mathcal{A}$). There exists a unique positive linear operator
$S$ in the commutant of the GNS representation corresponding to
$\varphi$ such that
\begin{equation}\label{eq3_8_1}
\ip{S\pi_{\varphi}(x)\xi_0}{\pi_{\varphi}(y)\xi_0}=\varphi'(y^*x),\quad(x,y\in\mathcal{A}).
\end{equation}
Conversely, for any positive operator $S$ in the commutant of
$\pi_{\varphi}(\mathcal{A})$, there is a unique positive
functional $\varphi'$ on $\mathcal{A}$ such that (\ref{eq3_8_1})
holds and $\varphi'\leq\|S\|\varphi$.
\end{corollary}

\begin{example}\label{ex3_9}{\bf [Groups]}
Take now $X=G$ a group. We know from theorem \ref{th1_7} that, if
$K:G\times G\rightarrow\mathbb{C}$ is positive definite and
satisfies
$$K(x,y)=K(zx,zy),\quad(x,y,z\in G),$$
then $K$ induces a unitary representation of $G$ on $H_K$. In the
next theorem we look at operators that intertwine these
representations.
\end{example}

\begin{theorem}\label{th3_10}
Suppose $G$ is a group and $K,K'$ are positive definite maps on
$G$ satisfying
$$K(x,y)=K(zx,zy),\,K'(x,y)=K'(zx,zy),\quad(x,y,z\in G).$$
Let $L:G\times G\rightarrow\mathbb{C}$ with $L^2\leq cKK'$ for
some $c>0$. If
\begin{equation}\label{eq3_10_1}
L(x,y)=L(zx,zy),\quad(x,y,z\in G),
\end{equation}
then there is a unique operator $S:H_K\rightarrow H_{K'}$ such
that
\begin{equation}\label{eq3_10_2}
S\pi_K(x)=\pi_{K'}(x)S,\quad(x\in G),
\end{equation}
\begin{equation}\label{eq3_10_3}
\ip{S\pi_K(x)\xi_0}{\pi_{K'}(y)\xi_0'}=L(x,y),\quad(x\in G).
\end{equation}
($(H_K,\pi_K,\xi_0),(H_{K'},\pi_{K'},\xi_0')$ are the unitary
representations of $G$ associated to $K$ and $K'$ respectively
(see theorem \ref{th1_7})). Moreover $\|S\|\leq\sqrt{c}$.
Conversely, if $S:H_K\rightarrow H_{K'}$ satisfies
(\ref{eq3_10_2}), then there is a unique $L$ that satisfies
(\ref{eq3_10_3}). In addition $L$ satisfies (\ref{eq3_10_1}) and
$L^2\leq\|S\|^2KK'$.
\end{theorem}
\begin{proof}
Recall that, if $v_K:G\rightarrow H_K$ and $v_{K'}:G\rightarrow
H_{K'}$ are the representations associated to $K$ by Kolmogorov's
theorem, then
$$\pi_K(x)\xi_0=v_K(x),\pi_{K'}(x)\xi_0'=v_{K'}(x),\quad(x\in G)$$
(see the proof of theorem \ref{th1_7}).
\par
Theorem \ref{th3_2} implies the existence of an operator
$S:H_K\rightarrow H_{K'}$ with $\|S\|\leq\sqrt{c}$ and
$$\ip{Sv_K(x)}{v_{K'}(y)}=L(x,y),\quad(x,y\in G).$$
The rest follows.
\end{proof}

\begin{corollary}\label{cor3_11}
Let $K,K'$ be two positive definite maps on the group $G$ that
satisfy
$$K(x,y)=K(zx,zy),K'(x,y)=K'(zx,zy),\quad(x,y,z\in G),$$
and $K'\leq cK$ for some $c>0$. Then there exists a unique
positive operator $S$ on $H_K$ in the commutant of the unitary
representation $\pi_K(G)$, such that
\begin{equation}\label{eq3_11_1}
\ip{S\pi_K(x)\xi_0}{\pi_K(y)\xi_0}=K'(x,y),\quad(x,y\in G).
\end{equation}
Conversely, for every positive operator $S$ in the commutant of
$\pi_K(G)$, there is a unique positive definite map $K'$ on $G$
that satisfies (\ref{eq3_11_1}) and
$$K'(x,y)=K'(zx,zy),\quad(x,y,z\in G).$$
\end{corollary}
\begin{proof}
It is an immediate conseqence of theorem \ref{th3_10}. It can also
be proved from corollary \ref{cor3_5}.
\end{proof}

\begin{example}\label{ex3_12}{\bf [Gabor type unitary systems]}
We proved in theorem \ref{th1_10} that, given a unimodular
$\lambda\in\mathbb{C}$ and a positive definite map $K$ on
$\mathbb{Z}^2$ that satisfies
\begin{equation}\label{eq3_12_1}
K((m+1,n),(m'+1,n'))=K((m,n),(m',n')),\quad(m,n,m',n'\in\mathbb{Z}),
\end{equation}
\begin{equation}\label{eq3_12_2}
K((m,n+1),(m',n'+1))=\lambda^{m-m'}K((m,n),(m',n')),\quad(m,n,m',n'\in\mathbb{Z}),
\end{equation}
there is a Gabor type unitary system on $H_K$ generated by two
unitaries $U_K$ and $V_K$. As the reader probably expects, we look
at the operators that intertwine these systems.
\end{example}

\begin{theorem}\label{th3_11}
Let $\lambda\in\mathbb{C}$, $|\lambda|=1$ and $K,K'$ positive
definite maps on $\mathbb{Z}^2$ satisfying the corresponding
relations (\ref{eq3_12_1}) and (\ref{eq3_12_2}). Let
$L:\mathbb{Z}^2\times\mathbb{Z}^2\rightarrow\mathbb{C}$ with the
property that $L^2\leq cKK'$ for some $c>0$. If $L$ satisfies the
relations (\ref{eq3_12_1}) and (\ref{eq3_12_2}), (with $K$
replaced by $L$, of course), then there is a unique operator
$S:H_K\rightarrow H_{K'}$ such that
\begin{equation}\label{eq3_13_1}
SU_K=U_{K'}S,\quad SV_K=V_{K'}S,
\end{equation}
\begin{equation}\label{eq3_13_2}
\ip{SU_K^mV_K^n\xi_0}{U_{K'}^{m'}V_{K'}^{n'}\xi_0'}=L((m,n),(m',n')),\quad(m,n,m',n'\in\mathbb{Z}).
\end{equation}
($(U_K,V_K,\xi_0),(U_{K'},V_{K'},\xi_0')$ are given by theorem
\ref{th1_10})). Moreover $\|S\|\leq\sqrt{c}$. Conversely, if
$S:H_K\rightarrow H_{K'}$ satisfies (\ref{eq3_13_1}), then there
exists a unique
$L:\mathbb{Z}^2\times\mathbb{Z}^2\rightarrow\mathbb{C}$ that
verifies (\ref{eq3_13_2}) and in addition $L$ will verify
(\ref{eq3_12_1}) and (\ref{eq3_12_2}), too and $L^2\leq\|S\|^2
KK'$.
\end{theorem}
\begin{proof}
We recall that $U_K^mV_K^n\xi_0=v_K(m,n)$ and similarly for $K'$,
$(m,n\in\mathbb{Z})$ (see the proof of theorem \ref{th1_10}).
Theorem \ref{th3_2} shows that there is an operator
$S:H_K\rightarrow H_{K'}$ with $\|S\|\leq\sqrt{c}$ and such that
(\ref{eq3_13_2}) holds. We need to check (\ref{eq3_13_1}). Take
$m,n,m',n'\in\mathbb{Z}$ and compute:
$$\ip{SU_KU_K^mV_K^n\xi_0}{U_{K'}^{m'}V_{K'}^{n'}\xi_0'}=L((m+1,n),(m',n'))$$
$$\ip{U_{K'}SU_K^mV_K^n\xi_0}{U_{K'}^{m'}V_{K'}^{n'}\xi_0'}=\ip{SU_K^mV_K^n\xi_0}{U_{K'}^{m'-1}V_{K'}^{n'}\xi_0'}$$
$$=L((m,n),(m'-1,n'))=L((m+1,n),(m',n')).$$
The density of the linear spans of
$\{U_K^mV_K^n\xi_0\,|\,m,n\in\mathbb{Z}\}$ and
$\{U_{K'}^{m'}V_{K'}^{n'}\xi_0'\,|\,m',n'\in\mathbb{Z}\}$ implies
$SU_K=U_{K'}S$. A similar calculation shows that $SV_K=V_{K'}S$.
\par
The converse follows from theorem \ref{th3_2}: if $L$ is defined
by (\ref{eq3_13_2}), the only thing that remains to be verified is
that $L$ satisfies (\ref{eq3_12_1}) and (\ref{eq3_12_2}), but this
is a consequence of (\ref{eq3_13_1}) and $U_KV_K=\lambda V_KU_K$,
$U_{K'}V_{K'}=\lambda V_{K'}U_{K'}$.
\end{proof}

\begin{corollary}\label{cor3_14}
If $K,K'$ are positive definite maps on $\mathbb{Z}^2$ satisfying
the relations (\ref{eq3_12_1}) and (\ref{eq3_12_2}) and $K'\leq
cK$ then there is a unique positive definite operator $S$ on $H_K$
that commutes with $U_K$ and $V_K$ and
\begin{equation}\label{eq3_14_1}
\ip{SU_K^mV_K^n\xi_0}{U_K^{m'}V_K^{n'}\xi_0}=K'((m,n),(m',n')),\quad(m,n,m',n'\in\mathbb{Z}).
\end{equation}
Conversely, if $S$ is a positive operator that commutes with $U_K$
and $V_K$ then $K'$ defined by (\ref{eq3_14_1}) satisfies
(\ref{eq3_12_1}) and (\ref{eq3_12_2}).
\end{corollary}
\begin{proof}
The proof follows the same lines as before.
\end{proof}

\begin{remark}
Theorem \ref{th3_2} gives us a general existence result for
intertwining operators. The next theorems answer the question what
conditions should be imposed on $L$ to obtain that its associated
operator $S$ intertwines the extra structure existent on $H_K$? We
saw that for $C^*$-algebras the necessary and sufficient condition
is that $L(x,y)=\varphi_0(y^*x)$ for some linear $\varphi_0$, for
groups we must have $L(x,y)=L(zx,zy)$ and for Gabor type unitary
systems, $L$ must satisfy the relations (\ref{eq3_12_1}) and
(\ref{eq3_12_2}).
\end{remark}

\begin{example}\label{ex3_15}{\bf [Intertwiners of wavelet representations]}
We mentioned in example \ref{ex1_12} and theorem \ref{th1_13} how
wavelet representations can be associated to positive functions
$h\in\lonet$ with $R_{m_0,m_0}h=h$. In \cite{Dut1} and \cite{Dut2}
we studied the operators that intertwine these representations. We
indicate now how these can be connected to Kolmogorov's theorem.
So we will recall the results from \cite{Dut1} and we sketch the
proof based on theorem \ref{th3_2}.
\par
Given $h$ as in theorem \ref{th1_13} call $\left(U_h,\pi_h , H_h ,
\varphi_h\right)$ the cyclic representation of $\mathfrak{A}_N$
associated to $h$. Also, define the transfer operator associated
to a pair $m_0,m_0'\in\linft$ by
$$R_{m_0,m_0'}f(z)=\frac{1}{N}\sum_{w^N=z}m_0(w)\overline{m_0 '}(w)f(w),\quad(z\in\mathbb{T},f\in\lonet).$$
\end{example}

\begin{theorem}\label{th3_16}
Let $m_0,m_0'\in\linft$ be non-singular and $h,h'\in\lonet$,
$h,h'\geq0$, $R_{m_0,m_0}(h)=h$, $R_{m_0',m_0'}(h')=h'$. Let
$\left(U,\pi,H,\varphi\right)$, $\left(U',\pi ',H',\varphi
'\right)$
be the cyclic representations corresponding to $h$ and $h'$ respectively. \\
If $h_0\in\lonet$, $R_{m_0,m_0'}\left(h_0\right)=h_0$ and $\left|
h_0\right| ^2\leq chh'$ for some $c>0$ then there exists a unique
operator $S:H\rightarrow H'$ such that
\begin{equation}\label{eq3_16_1}
SU=U'S,\quad S\pi(f)=\pi '(f)S,\quad (f\in\linft)
\end{equation}
\begin{equation}\label{eq3_16_2}
\ip{S\pi (f)\varphi }{\varphi
'}=\int_{\mathbb{T}}fh_0\,d\mu,\quad(f\in\linft).
\end{equation}
Moreover $\left\| S\right\|\leq\sqrt{c}$. Conversely, if $S$ is an
operator that satisfies (\ref{eq3_16_1}), then there is a unique
$h_0\in\lonet$ with $R_{m_0,m_0'}h_0=h_0$ such that
(\ref{eq3_16_2}) holds. Moreover, $|h_0|^2\leq\|S\|^2 hh'$.
\end{theorem}

\begin{proof}
Define $X$ as in the proof of theorem \ref{th1_13}. For all
$(f,n),(g,m)\in\ X$, we want to obtain
\begin{align*}
L((f,n),(g,m))&=\ip{SU^{-n}\pi(f)\varphi}{U'^{-m}\pi'(g)\varphi'}\\
&=\ip{SU^m\pi(f)\varphi}{U'^n\pi'(g)\varphi'}\\
&=\ip{S\pi(f(z^{N^m})m_0^{(m)}(z))\varphi}{\pi'(g(z^{N^n})m_0^{(n)}(z))\varphi'}\\
&=\int_{\mathbb{T}}f(z^{N^m})m_0^{(m)}(z)\overline{g}(z^{N^n})\overline{m_0'^{(n)}}(z)h_0\,d\mu.
\end{align*}
Keep the first and the last terms of the equality and this defines
$L$. $L$ will give rise to $S$ by theorem \ref{th3_2}. For the
details of the required computations, see \cite{Dut1}. The
converse, can be also obtained from theorem \ref{th3_2}, but here
 the generality of theorem \ref{th3_2} isn't really needed.
\end{proof}

\section{\label{Frames}Frames and dilations}
Recall that a set $\{x_i\,|\,i\in I\}$ of vectors in a Hilbert
space $H$ is called a frame if there are two constants $A,B>0$
such that
$$A\|f\|^2\leq\sum_{i\in I}|\ip{f}{x_i}|^2\leq B\|f\|^2,\quad(f\in H).$$
If $A=B=1$ the set $\{x_i\,|\, i\in I\}$ is called a normalized
tight frame.
\par
Frames have been used extensively in applied mathematics for
signal processing and data compression. They play a central role
in wavelet theory and the analysis of Gabor systems.
\par
In \cite{HL} the normalized tight frames are interpreted as
projections of orthonormal bases and it is proved there that Gabor
type normalized tight frames can be dilated to Gabor type
orthonormal bases, and normalized tight frames generated by groups
can be dilated to orthonormal bases generated by the same group
(see theorem 3.8 and 4.8 in \cite{HL}). We will revisit these
theorems and show that they are immediate consequences of a
general result which proves that any normalized tight frame can be
dilated to an orthonormal basis in such a way that the extra
structure that may exists is preserved under the dilation.
\par
We begin with a proposition that establishes what positive
definite maps give rise to normalized tight frames when
represented on a Hilbert space.
\begin{proposition}\label{prop4_1}
Let $K$ be a positive definite map on a set $X$. Then
$\{v_K(x)\,|\, x\in X\}$ is a normalized tight frame if and only
if for all $x_i\in X,\xi_i\in\mathbb{C}$, $(i\in\{1,...,n\})$:
\begin{equation}\label{eq4_1}
\sum_{i,j=1}^nK(x_i,x_j)\xi_i\overline{\xi}_j=\sum_{x\in
X}\left|\sum_{i=1}^nK(x_i,x)\xi_i\right|^2.
\end{equation}
\end{proposition}
\begin{proof}
If $\{v_K(x)\,|\, x\in X\}$ is a normalized tight frame then take
$f=\sum_{i=1}^n\xi_i v_K(x_i)$. The fact that
\begin{equation}\label{eq4_2}
\|f\|^2=\sum_{x\in X}|\ip{f}{v_K(x)}|^2,
\end{equation}
translates into (\ref{eq4_1}). For the converse, we only need to
verify (\ref{eq4_2}) for $f$ in a dense subset of $H_K$ (see
\cite{HeWe} lemma 1.10). Since the linear span of $\{v_K(x)\,|\,
x\in X\}$ is dense in $H_K$, we can take $f=\sum_{i=1}^n\xi_i
v_K(x_i)$ and (\ref{eq4_2}) follows from (\ref{eq4_1}).
\end{proof}

\begin{definition}\label{def4_2}
A positive definite map $K$ on $X$ is called a NTF if and only if
$\{v_K(x)\,|\, x\in X\}$ is a normalized tight frame for $H_K$.
\end{definition}

Before we prove our general result we note that, if $\delta
:X\times X\rightarrow\mathbb{C}$ is defined by
$$\delta(x,y)=\left\{\begin{array}{lcr}
1&,&\mbox{ if }x=y\\
0&,&\mbox{ otherwise }
\end{array}\right.$$
then $\delta$ is a positive definite map and $\{v_{\delta}(x)\,|\,
x\in X\}$ is an orthonormal basis for $H_{\delta}$.

\begin{proposition}\label{prop4_3}
If $K$ is a NTF positive definite map on $X$ then $K\leq\delta$.
\end{proposition}
\begin{proof}
By definition, $\{v_K(x)\,|\, x\in X\}$ is a normalized tight
frame for $H_K$. Then, by \cite{HL} proposition 1.1, there exists
a Hilbert space $H$ containing $H_K$ as a subspace and an
orthonormal basis $\{e(x)\,|\, x\in X\}$ such that, if $P$ is the
projection onto $H_K$, then $Pe(x)=v_K(x)$ for all $x\in X$.
\par
Now take $x_i\in X$, $\xi_i\in\mathbb{C}$, $(i\in\{1,...,n\}$.
\begin{align*}
\sum_{i,j=1}^nK(x_i,x_j)\xi_i\overline{\xi}_j&=\ip{\sum_{i=1}^n\xi_i v_K(x_i)}{\sum_{i=1}^n\xi_i v_K(x_i)}\\
&=\ip{P(\sum_{i=1}^n\xi_ie(x_i))}{P(\sum_{i=1}^n\xi_ie(x_i))}\\
&\leq\|P\|\,\|\sum_{i=1}^n\xi_i e(x_i)\|^2=\sum_{i=1}^n|\xi_i|^2\\
&=\sum_{i,j=1}^n\delta(x_i,x_j)\xi_i\overline{\xi}_j.
\end{align*}
Therefore $K\leq\delta$.
\end{proof}

\begin{theorem}\label{th4_4}
If $K,K'$ are NTF positive definite maps on a countable set $X$,
$K\leq cK'$ for some $c>0$, then there exists an isometry
$W:H_K\rightarrow H_{K'}$, $W$ is induced by $K$, that is
$$\ip{Wv_K(x)}{v_{K'}(y)}=K(x,y),\quad(x,y\in X),$$
the projection $P$ onto $WH_K$ is also induced by $K$, i.e.
$$\ip{Pv_{K'}(x)}{v_K'(y)}=K(x,y),\quad(x,y\in X),$$
and $Pv_{K'}(x)=Wv_K(x)$.
\end{theorem}
\begin{proof}
Since $K\leq cK'$, by corollary \ref{cor3_5}, there exists a
positive operator $S$ on $H_{K'}$ such that
$$\ip{Sv_{K'}(x)}{v_{K'}(y)}=K(x,y),\quad(x,y\in X).$$
Since $S$ is positive, it has a positive square root
$S^{\frac{1}{2}}$. Then
$$\ip{S^{\frac{1}{2}}v_{K'}(x)}{S^{\frac{1}{2}}v_{K'}(y)}=\ip{Sv_{K'}(x)}{v_{K'}(y)}=K(x,y),\quad(x,y\in X).$$
Take
$$H=\overline{\mbox{span}}\{S^{\frac{1}{2}}v_{K'}(x)\,|\, x\in X\}.$$
By the uniqueness part of Kolmogorov's theorem, there is a unitary
$W:H_K\rightarrow H$ such that
$$Wv_K(x)=S^{\frac{1}{2}}v_{K'}(x),\quad(x\in X).$$
But then $\{S^{\frac{1}{2}}v_{K'}(x)\,|\, x\in X\}$ is a
normalized tight frame for $H$. Also, we know that
$\{v_{K'}(x)\,|\, x\in X\}$ is a normalized tight frame for
$H_{K'}$. So $S^{\frac{1}{2}}:H_{K'}\rightarrow H$ maps a
normalized tight frame to a normalized tight frame, therefore it
must be a co-isometry (see \cite{HL} proposition 1.9). It follows
that $S^{\frac{1}{2}}\left(S^{\frac{1}{2}}\right)^*:H\rightarrow
H$ is the identity on $H$ so $S$ is the identity on $H$.
\par
We also know that $\mbox{range}(S^{\frac{1}{2}})=\mbox{range}(S)$.
This implies that $S(Sv)=Sv$ for all $v\in H_{K'}$ and, as
$S\geq0$, $S$ is the projection onto $H$. Consequently, we also
have $S=S^{\frac{1}{2}}$ and everything follows now by an easy
computation:
$$Sv_{K'}(x)=S^{\frac{1}{2}}v_{K'}(x)=Wv_{K}(x),\quad(x\in X).$$
$$\ip{Wv_{K}(x)}{v_{K'}(y)}=\ip{Sv_{K'}(x)}{v_{K'}(y)}=K(x,y),\quad(x,y\in X).$$
\end{proof}

\begin{remark}\label{rem4_5}
Theorem \ref{th4_4} can be used to construct dilation theorems for
unitary systems (the reader should have in mind the specific
examples of groups and Gabor type unitary systems).
\par
Recall some definitions from \cite{HL}. If $\mathcal{U}$ is a
countable set of unitaries on a Hilbert space $H$, then $\xi\in H$
is called a complete wandering vector (complete normalized tight
frame vector) if $\{U\xi\,|\,U\in\mathcal{U}\}$ is an orthonormal
basis (normalized tight frame) for $H$. A dilation theorem will
take the following form:
\par
If $\mathcal{U}$ is a unitary system on a Hilbert space $H$ that
has a complete normalized tight frame vector $\eta$, then there is
a Hilbert space $H_1$ that contains $H$ and a unitary system
$\mathcal{U}_1$ on $H_1$ such that $\mathcal{U}_1$ has a complete
wandering vector $\xi$ and if $P$ is the projection onto $H$ then
$P\xi=\eta$, $P$ commutes with $\mathcal{U}_1$ and $U_1\mapsto
U_1|_H$ is an isomorphism of $\mathcal{U}_1$ onto $\mathcal{U}$.
\par
The proof will be guided by the following steps:\\
1. Construct
$K:\mathcal{U}\times\mathcal{U}\rightarrow\mathbb{C}$,
$K(x,y)=\ip{x\eta}{y\eta}$; then $K$ is an NTF positive definite
map
and $H_K=H$, $v_K(x)=x\eta$, $(x\in\mathcal{U})$ and $\mathcal{U}$ is the extra structure $\mathcal{U}_K$ induced by $K$.\\
2. Verify that $\delta:\mathcal{U}\times\mathcal{U}\rightarrow\mathbb{C}$ satisfies the required compatibility conditions with $\mathcal{U}$.\\
3. Construct $H_\delta$, $v_\delta$ and the additional structure
$\mathcal{U}_\delta$ with cyclic vector $\xi_\delta$  which is
a complete wandering vector for $\mathcal{U}_\delta$.\\
4. Since $K\leq\delta$ (proposition \ref{prop4_3}), according to
theorem \ref{th4_4} there is an isometry $W:H\rightarrow H_\delta$
which is induced by $K$; the projection $P$ onto $WH$ is also
induced by $K$ and $P\xi_\delta=\eta$. As $K$ is compatible with
the structure $\mathcal{U}$, $W$ will intertwine $\mathcal{U}$ and
$\mathcal{U}_\delta$ and $P$ commutes with $\mathcal{U}_\delta$.
So $WH$ is invariant for $\mathcal{U}_\delta$ and
$WUW^{-1}=U_\delta|_H$ for all $U\in\mathcal{U}$ ($U_\delta$ is
the unitary in $\mathcal{U}_\delta$ that corresponds to $U$ in the
representation). 5. Identify $H$ with $WH$ and everything will
follow.
\end{remark}
We will use the guidelines of remark \ref{rem4_5} to show how one
can obtain the dilation theorems 3.8 and 4.8 from \cite{HL} for
groups and Gabor type unitary systems.

\begin{theorem}\label{th4_6}\cite{HL}
Suppose $\mathcal{U}$ is a unitary group on $H$ with a complete
normalized tight frame vector $\eta$. Then there is a Hilbert
space $H_1$ containing $H$ and a unitary group $\mathcal{U}_1$
such that $\mathcal{U}_1$ has a complete wandering vector $\xi$,
if $P$ is the projection onto $H$ then $P$ commutes with
$\mathcal{U}_1$, $P\xi=\eta$ and $U_1\mapsto U_1|_H$ is an
isomorphism of $\mathcal{U}_1$ onto $\mathcal{U}$. Consequently,
$PU_1\xi=U_1|_H\eta$ for all $U_1\in\mathcal{U}_1$ (that is the
normalized tight frame $\{U\eta\,|\, U\in\mathcal{U}\}$ can be
dilated to the orthonormal basis $\{U_1\xi\,|\,
U_1\in\mathcal{U}_1\}$).
\end{theorem}
\begin{proof}
Define $K:\mathcal{U}\times\mathcal{U}\rightarrow\mathbb{C}$,
$K(x,y)=\ip{x\eta}{y\eta}$ for $x,y\in\mathcal{U}$. It is clear
that $K$ is an NTF positive definite map with
\begin{equation}\label{eq4_6_1}
K(zx,zy)=K(x,y),\quad(x,y,z\in\mathcal{U})
\end{equation}
and $H_K=H$, $v_K(x)=x\eta$ and the representation $\pi_K$ given
by theorem \ref{th1_7} is $\pi_K(x)=x$ for $x\in\mathcal{U}$.
\par
It is also clear that $\delta$ satisfies a relation of type
(\ref{eq4_6_1}) so it is compatible with the group structure and
by theorem \ref{th1_7} it induces a cyclic representation
$(H_\delta,\pi_\delta,\xi_\delta)$ of $\mathcal{U}$ with
$\xi_\delta=v_\delta(1)$ a complete wandering vector.
\par
By proposition \ref{prop4_3}, $K\leq\delta$. By theorem
\ref{th4_4} there is an isometry $W:H\rightarrow H_\delta$ which
is induced by $K$, the projection $P$ onto $WH$ is also induced by
$K$ and $P\xi_\delta=\eta$. Then, by theorem \ref{th3_10}, $W$ is
intertwining that is
$$Wx=\pi(x)W,\quad(x\in\mathcal{U}),$$
and $P$ is in the commutant of $\pi_\delta(\mathcal{U})$. So $WH$
is invariant for all $\pi_\delta(x)$, $x\in\mathcal{U}$ and
$WxW^{-1}=\pi_\delta(x)$ for $x\in\mathcal{U}$.
\par
Identify $H$ with $WH$ and define
$\mathcal{U}_1=\pi_\delta(\mathcal{U})$, $\xi=\xi_\delta$ and
everything follows.
\end{proof}

\begin{theorem}\label{th4_6_2}\cite{HL}
Let $\mathcal{U}=\{U^mV^n\,|\,m,n\in\mathbb{Z}\}$ be a Gabor type
unitary system associated to $\lambda$ on a Hilbert space $H$.
Suppose $\mathcal{U}$ has a complete normalized tight frame vector
$\eta\in H$. Then there is a Gabor type unitary system
$\mathcal{U}_1(=\{U_1^mV_1^n\,m,n\in\mathbb{Z}\})$ associated to
$\lambda$ on a Hilbert space $H_1$ containing $H$, such that
$\mathcal{U}_1$ has a complete wandering vector $\xi$ and if $P$
is the projection onto $H$ then $P$ commutes with $U_1$ and $V_1$,
$P\xi=\eta$ and $U=U_1|_H$ $V=V_1|_H$.
\end{theorem}
\begin{proof}
The proof is analogous to the proof of theorem \ref{th4_6}, the
only difference is to verify that $\delta$ satisfies the
compatibility relations (\ref{eq3_12_1}) and (\ref{eq3_12_2}) and
this is trivial.
\end{proof}

\section{\label{A dilation}A dilation theorem for wavelets}
Let us recall the algorithm for the construction of compactly
supported wavelets. For details we refer the reader to
\cite{Dau92} for the scale $N=2$ and to \cite{BraJo97} for
arbitrary scale $N$.
\par
One starts with the low-pass filter $m_0\in\ltwot$ which is a
trigonometric polynomial that satisfies $m_0(1)=\sqrt{N}$ and the
quadrature mirror filter condition
\begin{equation}\label{eq5_0_1}
\frac{1}{N}\sum_{w^N=z}|m_0|^2(w)=1,\quad(z\in\mathbb{T}).
\end{equation}
Then define the scaling function $\varphi\in\ltwor$ by taking the
inverse Fourier transform of
\begin{equation}\label{eq5_0_2}
\widehat{\varphi}(x)=\prod_{k=1}^\infty\frac{m_0\left(\frac{x}{N^k}\right)}{\sqrt{N}},\quad(x\in\mathbb{R}).
\end{equation}
\par
To construct wavelets one needs the high-pass filters
$m_1,...,m_{N-1}\in\linft$ such that the matrix
\begin{equation}\label{eq5_0_3}
\frac{1}{\sqrt{N}}\left(\begin{array}{cccc}
m_0(z)&m_0(\rho z)&\hdots&m_0(\rho^{N-1}z)\\
m_1(z)&m_1(\rho z)&\hdots&m_1(\rho^{N-1}z)\\
\vdots&\vdots&\vdots&\vdots\\
m_{N-1}(z)&m_{N-1}(\rho z)&\hdots&m_{N-1}(\rho^{N-1}z)
\end{array}\right) \mbox{ is unitary for a.e }z\in\mathbb{T}.
\end{equation}
($\rho=e^{\frac{2\pi i}{N}}$).
\par
When $N=2$ a choice for $m_1$ can be
$$m_1(z)=z\overline{m}_0(-z)f(z^2),\quad(z\in\mathbb{T}),$$
where $|f(z)|=1$ on $\mathbb{T}$.

\par
The wavelets are defined as follows:
\begin{equation}\label{eq5_0_4}
\widehat{\psi}_i(x)=\frac{m_i\left(\frac{x}{N}\right)}{\sqrt{N}}\widehat{\varphi}\left(\frac{x}{N}\right),\quad(x\in\mathbb{R},i\in\{1,...,N-1\}),
\end{equation}
or, in terms of the wavelet representation,
\begin{equation}\label{eq5_0_5}
\psi_i=U^{-1}\pi(m_i)\varphi,\quad(i\in\{1,...,N-1\}).
\end{equation}
\par
It is known that, in order to achieve orthogonality, extra
conditions must be imposed on $m_0$. If $R_{m_0,m_0}$ has only one
continuous fixed point (up to a multiplicative constant), the set
$$\{U^mT^n\psi_i\,|\, m,n\in\mathbb{Z},i\in\{1,...,N-1\}\}$$
is an orthonormal basis for $\ltwor$.
\par
However, when this extra condition is not satisfied, one still
gets good properties, namely, the fact that the above set is a
normalized tight frame for $\ltwor$. In the sequel, we show how
one can dilate this normalized tight frame to an orthonormal basis
in such a way that the multiresolution structure is preserved so
that "wavelets" in a space bigger then $\ltwor$ are obtained.
\par
We begin with a proposition that explains the multiresolution
structure of the cyclic representations presented in example
\ref{ex1_12}.
\par In the sequel we define $m_0\in\linft$ to be non-singular if
the set $\{z\in\mathbb{T}\,|\, m_0(z)=0\}$ has zero measure and
$|m_0|$ is not constant 1 a.e.
\begin{proposition}\label{prop5_1}
Let $m_0\in\linft$ be non-singular, $h\in\lonet$, $h\geq 0$ and
$R_{m_0,m_0}h=h$. Let $(U_h,\pi_h,H_h,\varphi_h)$ be the cyclic
representation associated to $h$. Define $T_h=\pi_h(z)$,
$$V_0^h=\overline{\mbox{span}}\{T_h^k\varphi_h\,|\,k\in\mathbb{Z}\},$$
$$V_j^h=U_h^{-j}V_0^h,\quad(j\in\mathbb{Z}).$$
Then
\begin{equation}\label{eq5_1_1}
U_hT_hU_h^{-1}=T_h^N,
\end{equation}
\begin{equation}\label{eq5_1_2}
V_j^h\subset V_{j+1}^h,\quad(j\in\mathbb{Z}),
\end{equation}
\begin{equation}\label{eq5_1_3}
\overline{\cup_{j\in\mathbb{Z}}V_j^h}=H_h,
\end{equation}
\begin{equation}\label{eq5_1_4}
\cap_{j\in\mathbb{Z}}V_j^h=\{0\}
\end{equation}
Assume $h=1$. Then
\begin{equation}\label{eq5_1_5}
\{T_h^k\varphi_h\,|\, k\in\mathbb{Z}\}\mbox{ is an orthonormal
basis for }V_0^h.
\end{equation}
If $m_1,..,m_{N-1}$ satisfy (\ref{eq5_0_3}) and
\begin{equation}\label{eq5_1_6}
\psi_i^h=U_h^{-1}\pi_h(m_i)\varphi_h,\quad(i\in\{1,...,N-1\}),
\end{equation}
then
\begin{equation}\label{eq5_1_7}
\{T_h^k\psi_i^h\,|\,k\in\mathbb{Z},i\in\{1,...,N-1\}\}\mbox{ is an
orthonormal basis for }V_1^h\ominus V_0^h
\end{equation}
and
\begin{equation}\label{eq5_1_8}
\{U_h^mT_h^n\psi_i^h\,|\,m,n\in\mathbb{Z},i\in\{1,...,N-1\}\}\mbox{
is an orthonormal basis for }\ltwor.
\end{equation}
\end{proposition}

\begin{proof}
(\ref{eq5_1_1}) follows from
$U_h\pi_h(f(z))U_h^{-1}=\pi_h(f(z^N))$ with $f(z)=z$.
\par
If $f(z)=\sum_{k=-p}^pa_kz^k$ is a trigonometric polynomial then
$$\pi_h(f)\varphi_h=\sum_{k=-p}^pa_kT_k^h\varphi_h\in V_0^h.$$
Each $f\in\linft$ is the pointwise limit of a uniformly bounded
sequence of trigonometric polynomial, hence, by lemma 2.8 in
\cite{Dut1}, $\pi_h(f)\varphi_h\in V_0^h$. Then
$$U_h\pi_h(f)\varphi_h=\pi_h(f(z^N))U_h\varphi_h=\pi_h(f(z^N)m_0(z))\varphi_h\in V_0^h$$
so $V_{-1}^h\subset V_0^h$ and this implies (\ref{eq5_1_2}).
\par
Also $U_h^{-n}\pi_h(f)\varphi_h\in V_n^h$ for all $f\in\linft$ and
$n\in\mathbb{Z}$ and (\ref{eq5_1_3}) follows by density.
(\ref{eq5_1_4}) is proved in theorem 5.6 from \cite{Jor}.
\par
If $h=1$ then, for $k\in\mathbb{Z}$,
$$\ip{T_h^k\varphi_h}{\varphi_h}=\ip{\pi_h(z^k)\varphi_h}{\varphi_h}=\int_{\mathbb{T}}z^k\,d\mu=\delta_{k,0},$$
so (\ref{eq5_1_5}) is valid.
\par
It remains to prove (\ref{eq5_1_7}) because (\ref{eq5_1_8})
follows from this immediately. The argument is essentially the one
in \cite{BraJo97} theorem 10.1. We will include it here to make
sure everything works.
\par
For $k,l\in\mathbb{Z}$ and $i,j\in\{1,...,N-1\}$ we have
\begin{align*}
\ip{T_h^k\psi_i^h}{T_h^l\psi_j^h}&=\ip{\pi_h(z^k)U_h^{-1}\pi_h(m_i)\varphi_h}{\pi_h(z^l)U_h^{-1}\pi_h(m_j)\varphi_h}\\
&=\ip{U_h^{-1}\pi_h(z^{Nk}m_i(z))\varphi_h}{U_h^{-1}\pi_h(z^{Nl}m_j(z))\varphi_h}\\
&=\int_{\mathbb{T}}z^{N(k-l)}m_i(z)\overline{m}_j(z)\,d\mu\\
&=\int_{\mathbb{T}}z^{k-l}\frac{1}{N}\sum_{w^N=z}m_i(w)\overline{m}_j(w)\,d\mu=\delta_{i,j}\delta_{k,l},
\end{align*}
for the last equality we used (\ref{eq5_0_3}). So
$$\{T_h^k\psi_i^h\,|\,k\in\mathbb{Z},i\in\{1,...,N-1\}\}$$
is an orthonormal set.
\par
Take, $U_h^{-1}\pi_h(m)\varphi_h\in V_1^h$, $m\in\linft$ (vectors
of this form are dense in $V_1^h$).
\\ $U_h^{-1}\pi_h(m)\varphi_h\perp V_0^h$ is equivalent to,
for all $f\in\linft$:
\begin{align*}
0&=\ip{U_h^{-1}\pi_h(m)\varphi_h}{\pi_h(f)\varphi_h}=\ip{\pi_h(m)\varphi_h}{\pi_h(f(z^N)m_0(z))\varphi_h}\\
&=\int_{\mathbb{T}}m(z)\overline{f}(z^N)\overline{m}_0(z)\,d\mu\\
&=\int_{\mathbb{T}}\frac{1}{N}\sum_{w^N=z}m(w)\overline{m}_0(w)\overline{f}(z)\,d\mu
\end{align*}
which is equivalent to
$$\frac{1}{N}\sum_{w^N=z}m(w)\overline{m}_0(w)=0\mbox{ a.e. on }\mathbb{T}.$$
This shows in particular that $\psi_i^h\perp V_0^h$ for all
$i\in\{1,...,N-1\}$. Also, the vector
$$\vec{m}(z)=(m(z),m(\rho z),...,m(\rho^{N-1}z))$$
($\rho=e^{-2\pi i/N}$) must be perpendicular to the vector
$$\vec{m}_0(z)=(m_0(z),m_0(\rho z),...,m_0(\rho^{N-1}z))$$
for almost all $z$, so
$$\vec{m}(z)=\sum_{k=1}^{N-1}\mu_k(z)\vec{m}_k(z)$$
where $\mu_k(z)=\ip{\vec{m}(z)}{\vec{m}_k(z)}$ ( which shows that
$\mu_k\in\linft$).
\par
Since $\vec{m}(\rho z)$ is a circular permutation of $\vec{m}(z)$,
it follows that we must have $\mu_k(\rho z)=\mu_k(z)$, that is
$\mu_k(z)=\lambda_k(z^N)$ for some $\lambda_k\in\linft$.
\par
Then
$$m(z)=\sum_{k=1}^{N-1}\lambda_k(z^N)m_k(z),\quad(z\in\mathbb{T}),$$
and we compute
$$U_h^{-1}\pi_h(m)\varphi_h=U_h^{-1}\pi_h(\sum_{k=1}^{N-1}\lambda_k(z^N)m_k(z))\varphi_h=\sum_{k=1}^{N-1}\pi_h(\lambda_k)\psi_k^h.$$
and this shows that
$$U_h^{-1}\pi_h(m)\varphi_h\in\overline{\mbox{span}}\{T_h^k\psi_i^h\,|\,k\in\mathbb{Z},i\in\{1,...,N-1\}\}$$
by an argument similar to the one used in the begining of the
proof (now for $\psi_i^h$ instead of $\varphi_h$). This completes
the proof of (\ref{eq5_1_7}).
\end{proof}

\par
Motivated by the discussion in the begining of this section, we
give a dilation theorem for wavelets. The theorem describes how
one can dilate a normalized tight frame wavelet to an orthonormal
wavelet in a bigger space.

\begin{theorem}\label{th5_2}.
Let $m_0\in\linft$ be a non-singular filter with $R_{m_0,m_0}1=1$,
$h\in\linft$, $h\geq0$, $R_{m_0,m_0}h=h$, and consider
$(U_h,\pi_h,H_h,\varphi_h)$, the cyclic representation associated
to $h$. Assume also there are given filters
$m_1,...,m_{N-1}\in\linft$ such that (\ref{eq5_0_3}) holds, define
$\psi_i^h$ as in (\ref{eq5_1_6}) ($i\in\{1,...,N-1\}$) and suppose
$$\{U_h^mT_h^n\psi_i^h\,|\,m,n\in\mathbb{Z},i\in\{1,...,N-1\}\}$$
is a normalized tight frame for $H_h$ ($T_h=\pi_h(z)$).
\par
Then, if $(U_1,\pi_1,H_1,\varphi_1)$ is the cyclic representation
associated to the constant function $1$, then there exists an
isometry $W:H_h\rightarrow H_1$ with the following properties:
\begin{enumerate}
\item $WU_h=U_1W$, $W\pi_h(f)=\pi_1(f)W$ for all $f\in\linft$;
\item If $P$ is the projection onto $WH_h$ then
\begin{equation}\label{eq5_2_1}
PU_1=U_1P,\quad P\pi_1(f)=\pi_1(f)P,\quad(f\in\linft);
\end{equation}
\begin{equation}\label{eq5_2_2}
P\varphi_1=W\varphi_h.
\end{equation}
\item If $\psi_i^1=U_1^{-1}\pi_1(m_i)\varphi_1$,
$i\in\{1,...,N-1\}$ then
\begin{equation}\label{eq5_2_3}
\{U_1^mT_1^n\psi_i^1\,|\,m,n\in\mathbb{Z},i\in\{1,...,N-1\}\}\mbox{
is an orthonormal basis for }H_1,
\end{equation}
where $T_1=\pi_1(z)$;
\begin{equation}\label{eq5_2_4}
P\psi_i^1=W\psi_i^h,\quad(i\in\{1,...,N-1\}.
\end{equation}
\end{enumerate}
\end{theorem}

\begin{proof}
The proof is similar to the one of theorem \ref{th4_4} but some
additional arguments are needed. Since $h\in\linft$, we have
$|h|^2\leq \|h\|_\infty 1$ so, by theorem \ref{th3_16} there is a
positive operator $S$ on $H_1$ that commutes with $U_1$ and
$\pi_1$ and
$$\ip{S\pi_1(f)\varphi_1}{\varphi_1}=\int_{\mathbb{T}}fh\,d\mu,\quad(f\in\linft).$$
\par
$S$ has a positive sqare root $S^{\frac12}$ that commutes with
$U_1$ and $\pi_1$. Also the projection $P$ onto the range $H$ of
$S^{\frac12}$ must commute with $U_1$ and $\pi_1$. Then we can
restrict $\pi_1$ and $U_1$ to $H$ and
$$\ip{\pi_1(f)S^{\frac12}\varphi_1}{S^{\frac12}\varphi_1}=\int_{\mathbb{T}}fh\,d\mu,$$
$$U_1S^{\frac12}\varphi_1=\pi_1(m_0)S^{\frac12}\varphi_1.$$
The uniqueness part of theorem \ref{th1_13} implies that there is
a unitary $W$ from $H_h$ to $H$ with
$W\varphi_h=S^{\frac12}\varphi_1$, $WU_h=U_1W$ and
$W\pi_h(f)=\pi_1(f)W$ for $f\in\linft$. From these commuting
properties of $W$ and $P$ it follows that
$$W(U_h^mT_h^n\psi_i^h)=S^{\frac12}(U_1^mT_1^n\psi_1^h),\quad(m,n\in\mathbb{Z},i\in\{1,...,N-1\}).$$
Hence, $S^{\frac12}$ maps an orthonormal basis to a normalized
tight frame so it must be a co-isometry. Then, proceeding as in
the proof of theorem \ref{th4_4} we get that $S=S^{\frac12}=P$ and
everything follows.
\end{proof}

\par
When $m_0$ is a regular filter (we will give the precise meaning
of that in a moment), we can really get our hands on the abstract
cyclic representation associated to the constant function $1$ so
that we obtain a very concrete dilation theorem for non-orthogonal
wavelets in $\ltwor$. The construction is given in the next
theorem and it is based on the results presented in \cite{Dut2}.
\par
Before we state the result, some definitions are needed. A vector
$(z_1,z_2,...,z_p)$ is called an $m_0$-cycle if
$z_1^N=z_2,z_2^N=z_3,...,z_p^N=z_1$, $z_i$ are distinct and
$|m_0(z_i)|=\sqrt{N}$ for all $i\in\{1,...,p\}$.
\par
For $f\in\linft$ and $z_0\in\mathbb{T}$ define
$\alpha_{z_0}(f)(z)=f(zz_0)$ for $z\in\mathbb{T}$. For
$n\in\mathbb{N}$
$$m_0^{(n)}(z)=m_0(z)m_0(z^N)...m_0(z^{N^{n-1}}),\quad(z\in\mathbb{T}).$$

\begin{theorem}\label{th5_3}
Let $m_0$ be a Lipschitz function on $\mathbb{T}$ with finitely
many zeroes $m_0(1)=\sqrt{N}$, $R_{m_0,m_0}1=1$. Let
$C_j=(z_{1,j},...,z_{p_j,j})$ be the $m_0$-cycles,
$j\in\{1,...,n\}$, $m_0(z_{k,j})=\sqrt{N}e^{i\theta_{k,j}}$ for
all $k\in\{1,...,p_j\}$, $j\in\{1,...,n\}$,
$\theta_j=\theta_{1,j}+...+\theta_{p_j,j}$.
\par
For each $j\in\{1,...,n\}$ define: $H_j=\ltwor^{p_j}$,
$U_j:H_j\rightarrow H_j$
$$U_j(\xi_1,...,\xi_{p_j})=\left(e^{i\theta_{1,j}}U\xi_2,...,e^{i\theta_{p_j-1,j}}U\xi_{p_j},e^{i\theta_{p_j,j}}U\xi_1\right),$$
where $U\xi(x)=\frac{1}{\sqrt{N}}f\left(\frac{x}{N}\right)$ for
$\xi\in\ltwor$.
\par
For $f\in\linft$
$$\pi_j(f)(\xi_1,...,\xi_{p_j})=\left(\pi\left(\alpha_{z_{1,j}}(f)\right)\left(\xi_1\right),...,\pi\left(\alpha_{z_{p_j,j}}(f)\right)\left(\xi_p\right)\right),$$
where $\pi$ is the representation on $\ltwor$ defined in example
\ref{ex1_12}.
$$\widehat{\varphi}_{k,j}(x)=\prod_{l=1}^{\infty}\frac{e^{-i\theta_j}\alpha_{z_{k,j}}\left(m_0^{(p_j)}\right)\left(\frac{x}{N^{lp_j}}\right)}{\sqrt{N^{p_j}}},\quad (l\in\{1,...,p_j\}),$$
$$\varphi_j=(\varphi_{1,j},...,\varphi_{p_j,j}).$$
Finally, define $H_0=H_1\oplus...\oplus H_n$,
$U_0=U_1\oplus...\oplus U_n$, $\pi_0(f)=\pi_1(f)\oplus...\pi_n(f)$
for $f\in\linft$ and
$\varphi_0=\varphi_1\oplus...\oplus\varphi_n$. Then
$(U_0,\pi_0,H_0,\varphi_0)$ is the cyclic representation
associated to the constant function $1$.
\par
Also if $C_1$ is the trivial $m_0$-cycle $C_1=(1)$ then
$H_1=\ltwor$, $U_1=U$ ,$\pi_1=\pi$ and
$$\widehat{\varphi}_1(x)=\prod_{l=1}^\infty\frac{m_0\left(\frac{x}{N^l}\right)}{\sqrt{N}},\quad(x\in\mathbb{R}),$$
so $(U_1,\pi_1,H_1,\varphi_1)$ is the usual wavelet representation
on $\ltwor$.
\par
If $T_j=\pi_j(z)$, for $j\in\{1,...,n\}$ then
$$T_j(\xi_1,...,\xi_{p_j})(x)=(z_{1,j}T\xi_1,...,z_{p_j,j}T\xi_{p_j}),$$
where $T\xi(x)=\xi(x-1)$ for $\xi\in\ltwor$.
\par
Assume $m_1,...,m_{N-1}\in\linft$ satisfy (\ref{eq5_0_3}). Define
$$\psi_i^1=U_1^{-1}\pi_1(m_i)\varphi_1(\in\ltwor),\,\psi_i^0=U_0^{-1}\pi_0(m_i)\varphi_0,\,i\in\{1,...,N-1\}$$
and let $P_1$ be the projection from $H_0$ onto $H_1$ and
$T_0=T_1\oplus...\oplus T_n$. Then
\begin{equation}\label{eq5_3_1}
P_1U_0=U_0P_1,\quad P_1T_0=T_0P_1,\quad
P_1\pi_0(f)=\pi_0(f)P_1,\quad(f\in\linft);
\end{equation}
\begin{equation}\label{eq5_3_2}
U_0|_{H_1}=U_1(=U),\,
T_0|{H_1}=T_1(=T),\,\pi_0(f)|_{H_1}=\pi_1(f)(=\pi(f)),\quad(f\in\linft);
\end{equation}
\begin{equation}\label{eq5_3_3}
P_1\varphi_0=\varphi_1,\quad
P_1\psi_i^0=\psi_i^1,\quad(i\in\{1,...,N-1\});
\end{equation}
\begin{equation}\label{eq5_3_3_1}
U_0\varphi_0=\pi_0(m_0)\varphi_0,\quad
U_1\varphi_1=\pi_1(m_0)\varphi_1;
\end{equation}
\begin{equation}\label{eq5_3_3_2}
\{T_0^k\varphi_0\,|\,k\in\mathbb{Z}\}\mbox{ is an orthonormal
set;}
\end{equation}
\begin{equation}\label{eq5_3_4}
\{U_0^mT_0^n\psi_i^0\,|\,
m,n\in\mathbb{Z},i\in\{1,...,N-1\}\}\mbox{ is an orthonormal basis
for }H_0;
\end{equation}
\begin{equation}\label{eq5_3_5}
\{U_1^mT_1^n\psi_i^1\,|\,
m,n\in\mathbb{Z},i\in\{1,...,N-1\}\}\mbox{ is a normalized tight }
\end{equation}
$$\mbox{frame for }H_1=\ltwor.$$
\end{theorem}

\begin{proof}
The fact that the cyclic representation associated to the constant
function $1$ is proved in \cite{Dut2}, one needs only to take the
inverse Fourier transform of the representation presented there to
obtain the one described here. Then (\ref{eq5_3_1}) and
(\ref{eq5_3_2}) follow trivially from the definition,
(\ref{eq5_3_3}) follows from the definition and the commuting
properties of $P_1$, (\ref{eq5_3_3_1}) is included in the
definition of the cyclic representation, (\ref{eq5_3_3_2}) and
(\ref{eq5_3_4}) are consequences of proposition \ref{prop5_1} and
(\ref{eq5_3_5}) (which is also well known, see \cite{Dau92} or
\cite{BraJo97}) follows from the fact that the projection of an
orthonormal basis is a normalized tight frame (see \cite{HL}).
\end{proof}

\begin{example}\label{ex5_4}
We apply theorem \ref{th5_3} to the low-pass filter
$$m_0(z)=\frac{1+z^3}{\sqrt{2}}=\sqrt{2}e^{-\frac{3i\theta}{2}}\cos\left(\frac{3\theta}{2}\right),\quad(z=e^{-i\theta}\in\mathbb{T})$$
which is known to give non-orthogonal wavelets. The scale $N=2$.
Some short computations show that $m_0(1)=\sqrt{2}$,
$R_{m_0,m_0}1=1$. The $m_0$-cycles are
$$C_1=(z_{1,1}=1), C_2=(z_{2,1}=e^{2\pi i/3},z_{2,2}=e^{4\pi i/3})$$
$p_1=1,p_2=2$, $m_0(z_{1,1})=m_0(z_{2,1})=m_0(z_{2,2})=\sqrt{2}$
so $\theta_{1,1}=\theta_{2,1}=\theta_{2,2}=0$ and
$\theta_1=\theta_2=0$.
$$U_0:\ltwor^3\rightarrow\ltwor^3,\, U_0(\xi_1,\xi_2,\xi_3)=(U\xi_1,U\xi_3,U\xi_2),$$
$$T_0:\ltwor^3\rightarrow\ltwor^3,\, T_0(\xi_1,\xi_2,\xi_3)=(T\xi_1,e^{2\pi i/3}T\xi_2,e^{4\pi i/3}T\xi_3).$$
Then, as $\alpha_{z_{1,1}}(m_0)=m_0$,
$\alpha_{z_{2,1}}(m_0)=\alpha_{z_{2,2}}(m_0)=m_0$,
$$\widehat{\varphi}_{1,1}(x)=\prod_{l=1}^\infty e^{-\frac{3ix}{2^{l+1}}}\cos\left(\frac{3x}{2^{l+1}}\right)=e^{-\frac{3ix}{2}}\frac{\sin\left(\frac{3x}{2}\right)}{\frac{3x}{2}},\quad(x\in\mathbb{R}),$$
$$
\widehat{\varphi}_{2,1}(x)=\prod_{l=1}^\infty\frac{m_0^{(2)}\left(\frac{x}{2^{2l}}\right)}{\sqrt{2^2}}
=\prod_{l=1}^\infty\frac{m_0\left(\frac{x}{2^{2l}}\right)m_0\left(\frac{x}{2^{2l-1}}\right)}{\sqrt{2^2}}
=\widehat{\varphi}_{1,1}(x),
$$
and similarly for $\widehat{\varphi}_{2,2}$. Hence,
$\varphi_{1,1}=\varphi_{2,1}=\varphi_{2,2}=:\varphi=\frac{1}{3}\chi_{[0,3)}$,
and $\varphi_0=(\varphi,\varphi,\varphi)$.
\par
To construct the wavelet we can pick
$$m_1(z)=\frac{1-z^3}{\sqrt{2}},\quad(z\in\mathbb{T}).$$
Then the wavelet $\psi_0=(\psi_1,\psi_2,\psi_3)$ is given by
$$U_0\psi_0=\frac{1}{\sqrt{2}}(\varphi_0-T_0^3\varphi_0)$$
so
$$\psi_1=\psi_2=\psi_3=:\psi=\frac{1}{3}(\chi_{[0,\frac{3}{2})}-\chi_{[\frac{3}{2},1)}).$$
and
$$\{U_0^mT_0^n\psi_0\,|\,m,n\in\mathbb{Z}\}$$
is an orthonormal basis for $\ltwor^3$ which dilates the
normalized tight frame of $\ltwor$
$$\{U^mT^n\psi\,|\,m,n\in\mathbb{Z}\}.$$
\end{example}

\begin{acknowledgements}
The author wants to thank professor \c Serban Str\u atil\u a for
pointing out the connection between the GNS construction and
Kolmogorov's theorem. This was the starting point and the key idea
of this paper. Also many thanks to professor Palle Jorgensen for
his suggestions and his constant support.
\end{acknowledgements}

\end{document}